\documentclass[nofootinbib,superscriptaddress]{revtex4}
\usepackage{graphicx}
\usepackage{amsmath,amsthm,amssymb}
\begin{document}
\title{Self-similar solutions of semilinear wave equations\\ with a focusing nonlinearity}

\author{Piotr Bizo\'n}
\affiliation{Institute of Physics,
   Jagiellonian University, Krak\'ow, Poland}
\author{Dieter Maison}
\affiliation{Max-Planck-Institute f\"ur Physik,
Werner-Heisenberg-Institut,
   Munich, Germany}
\author{Arthur Wasserman} \affiliation{Department of Mathematics,
   University of Michigan, Ann Arbor, MI}

\begin{abstract}
\vskip 0.3cm \noindent We prove that in three space dimensions a
nonlinear wave equation $u_{tt}-\Delta u = u^p$ with $p\geq 7$ being
an odd integer has a countable family of regular spherically
symmetric self-similar solutions.
\end{abstract}

\maketitle

\subsection*{1. Introduction}
 An important feature of many nonlinear wave equations is that
their solutions, corresponding to smooth initial data, may form
singularities after a finite time. Such a phenomenon, usually
referred to as blowup, has been a subject of intensive studies
beginning with the pioneering works by Keller \cite{k}, John
\cite{j} and Glassey \cite{g} (we refer an interested reader to the
excellent online review of this subject with the complete
bibliography \cite{dw}).

 In this paper we  consider the semilinear wave equation
with a power nonlinearity
\begin{equation}\label{deq} \Phi_{tt}-\Delta\Phi-\Phi^p=0, \quad
\Phi=\Phi(t,x),\quad x\in R^3, \end{equation}
where $p\geq 7$ is an odd integer. The sign of the nonlinear term
corresponds to focusing which means that it tends to magnify the
amplitude of the wave. If $\Phi$ is small this term is negligible
and the evolution is essentially linear leading to dispersion.
However, if $\Phi$ is large the dispersive effect of the laplacian
may be overcome by the focusing effect of the nonlinearity and a
singularity can form. Actually, neglecting the laplacian altogether
 and solving the ordinary differential equation $\Phi_{tt}=\Phi^p$
 one gets the exact, homogeneous in space, solution
\begin{equation}\label{phi0}
    \Phi_0=\frac{b_0}{(T-t)^{\alpha}}, \quad
    b_0=\left[\frac{2(p+1)}{(p-1)^2}\right]^{\frac{1}{p-1}}, \quad
    \alpha=\frac{2}{p-1},\quad T>0,
\end{equation}
which blows up as $t\rightarrow T$. By the finite speed of
propagation one can truncate this solution in space to get a
solution with compactly supported initial data which blows up in
finite time. There is theoretical \cite{kl} and numerical \cite{bct}
evidence that the solution $\Phi_0$ determines the leading order
asymptotics of blowup for generic large initial data. However, for
specially prepared initial data, in particular for data fine-tuned
to the threshold for blowup, singularities may have a different form
which is given by self-similar solutions of equation (\ref{deq}).
Such solutions were found numerically in \cite{bct}. The aim of this
paper is to give a rigorous proof of their existence and discuss
their properties.

By definition, self-similar solutions are invariant under rescaling
\begin{equation}\label{css}    \Phi(t,x)\rightarrow
\Phi_{\lambda}(t,x)=\lambda^{-\alpha} \Phi(t/\lambda,x/\lambda),
\end{equation}
hence in the spherically symmetric case they have the
form
\begin{equation}\label{ansatz}
\Phi(t,r)=(T-t)^{-\alpha} u(\rho), \quad \rho=\frac{r}{T-t},
\end{equation}
where $T$ is a positive constant. Note that each  self-similar
solution, if it is regular for $t<T$, provides an explicit example
of regular initial data developing a singularity  in a finite time.
Substituting the ansatz (\ref{ansatz})  into equation (\ref{deq}) we
obtain  the ordinary differential equation for the similarity
profile $u(\rho)$
\begin{equation}\label{sseq}
(1-\rho^2)u''+\left(\frac{2}{\rho}-(2+2\alpha)\rho\right)u'-\alpha(\alpha+1)
u+u^p=0\,.
\end{equation}
It is easy to see that this equation has the constant solution
$u_0(\rho)=b_0$ which, of course, corresponds to the homogeneous
solution $\Phi_0$ of equation (\ref{deq}). As mentioned above,
numerical results indicate that, besides $u_0$, there exist also
nontrivial regular solutions of equation (\ref{sseq}). We remark
that equation (\ref{sseq}) has been studied by Kavian and Weissler
\cite{kw} who made many interesting observations about the behavior
of solutions but unfortunately they imposed very restrictive
fall-off conditions at infinity which excluded nontrivial solutions.
However, as long as one wants to have an example of blowup, due to
the finite speed of propagation, only the behavior of solutions
inside the past light cone of the blowup point $(t=T,r=0)$ is
relevant, which corresponds to  the interval $0\leq \rho \leq 1$.
 As we
shall see below regular self-similar solutions are parametrized by
their values at the endpoints $c=u(0)$ and $b=u(1)$. Due to the
singular nature of equation (\ref{sseq}) at $\rho=0$ and $\rho=1$,
the generic solution obeying the regularity  condition at $\rho=0$
becomes singular at $\rho=1$ and vice versa, hence globally regular
solutions can exist only for discrete values of the parameters $b$
and $c$. The goal of this paper is to show that there exists an
infinite sequence of pairs $(b_n,c_n)$ which give rise to globally
regular solutions. To this end, in section~2 we first prove
existence of local solutions near the endpoints. Then, in section~3,
we show, using a shooting argument, that for a discrete set of
values of the parameters these local solutions match smoothly at a
midpoint. In section~4 we show that the solutions constructed in
section~3 can be extended  beyond the past light cone to infinity.
Finally, in section~5 we derive some scaling properties of the
shooting parameters. In the appendix we present an alternative
rigorous proof of existence.
\subsection*{2. Local existence}
\noindent In the first step we will analyze the behavior of
solutions near the boundary points $\rho=0$ and $\rho=1$ and prove
the local existence of one-parameter families of regular solutions
there.

Near $\rho=0$ we introduce a new variable $v=u'$ and rewrite
equation (\ref{sseq}) as the first order system
\begin{subequations}
\begin{eqnarray}
\rho u'&=&\rho v,\\
\rho v'&=&-2v+\frac{\rho}{1-\rho^2}\Bigl(\frac{4}{p-1}\rho v+
     \frac{2(p+1)}{(p-1)^2}u-u^p\Bigr)\,.
\end{eqnarray}
\end{subequations}
This system has the form required by Proposition~1 of \cite{bfm}
guaranteeing the existence of a one-parameter family of regular
solutions near $\rho=0$ such that
\begin{equation}\label{local0}
    u(\rho)=c+\frac{1}{3} \left(\frac{p+1}{(p-1)^2} c - \frac{1}{2}
    c^p\right) \rho^2 + O(\rho^4).
\end{equation}
Similarly, near $\rho=1$ we put $s=1-\rho, \bar u(s)=u(\rho),\bar
v(s)=\bar u'(s)$ to get
\begin{subequations}
\begin{eqnarray}\label{rho1}
s \bar u'&=&s \bar v,\\
s \bar v'&=&\frac{1}{(2-s)}\Biggl(\Bigl(\frac{2}{1-s}
         -\frac{2(p+1)(1-s)}{p-1}\Bigr)\bar v+
             \frac{2(p+1)}{(p-1)^2}\bar u-{\bar u}^p\Biggr)\,.
\end{eqnarray}
\end{subequations}
In order to bring this into the form required by Proposition~1 of
\cite{bfm} we introduce $\tilde v=\bar v-\frac{p+1}{2(p-1)}\bar
u+\frac{p-1}{4}{\bar u}^p$ and rewrite the system (8) as
\begin{subequations}
\begin{eqnarray}
s\bar  u'&=&s\Bigl(\tilde v+\frac{p+1}{2(p-1)}\bar u-\frac{p-1}{4}{\bar u}^p\Bigr),\\
s\tilde v' &=&-\frac{2}{p-1}\tilde v+sf(\bar u,\tilde v,s),
\end{eqnarray}
\end{subequations}
where the function $f(\bar u,\tilde v,s)$ is analytic in $\bar u$,
$\tilde v$ and $s$ near $s=0$.
We conclude that there is a one-parameter family of regular
solutions near $\rho=1$ such that
\begin{equation}\label{local1}
    u(\rho)=b + \frac{1}{2} \left(\frac{1}{2} (p-1) b^p -
    \frac{p+1}{p-1} b \right) (\rho-1) + O((\rho-1)^2).
\end{equation}
 The singular solutions behave like
$u\sim (1-\rho)^{\frac{p-3}{p-1}}$ near $\rho=1$ leading to a
singularity of $u'$. For later use it is convenient to introduce the
variables $\sigma=(1-\rho)^{2/(p-1)},\psi(\sigma)=u(\rho)$ and
$\theta(\sigma)=\sigma u'$. The system (8) may then be written in
the form
\begin{subequations}
\begin{eqnarray}\label{rho1r}
\psi'&=&-\frac{p-1}{2}\sigma^{\frac{p-5}{2}} \theta,\\
\theta' &=&g(\psi,\theta,\sigma), \end{eqnarray}
\end{subequations}
where the function $g(\psi,\theta,\sigma)$ is analytic in $\psi,
\theta$ and $\sigma$ near $\sigma=0$. This implies that there is a
two-parameter family of regular solutions parametrized by the values
$b=\psi(0)=u(1)$ and $d=\theta(0)$. These solutions are analytic as
functions of $b$, $d$ and $\sigma$ in a neighborhood of $\sigma=0$.
Thus, the general solution is in fact regular at $\rho=1$ as a
function of $\sigma=(1-\rho)^{2/(p-1)}$ with a convergent Taylor
series at $\sigma=0$.
\subsection*{3. Global behavior and proof of existence}
\noindent The main result of this paper is: \vskip 0.2cm \noindent
\textbf{Theorem~1.} For any odd $p\geq 7$ and any non-negative
integer $n$, equation (\ref{sseq}) on the interval $0\leq \rho\leq
1$ has an analytic solution $u_n(\rho)$ which satisfies the boundary
conditions (\ref{local0}) and (\ref{local1}), and for which the
function $w_n=u_n/u_{\infty}-1$, with $u_{\infty}$ defined by
(\ref{uinf}), has exactly $n+1$ zeros. \vskip 0.2cm \noindent
\emph{Remark:} In this section we give a "physicist" proof (no
epsilon-delta stuff) based on scaling arguments. We believe that
sacrificing slightly the mathematical rigour we gain better
understandability of the basic mechanism which is responsible for
the existence and structure of solutions. An alternative rigorous
proof is presented in the appendix. \vskip 0.2cm
 \noindent \emph{Proof:} First, we need to show that
solutions which are regular at $\rho=0$ or $\rho=1$ remain bounded
on the whole interval $0<\rho<1$. This is conveniently proven using
a suitable Lyapunov function \cite{H} defined as
\begin{equation}\label{liap}
H=(1-\rho^2)\frac{u'^2}{2}+\frac{u^{p+1}}{p+1}-\frac{(p+1)u^2}{(p-1)^2}\,,
\end{equation}
obeying
\begin{equation}
H(\rho)\geq
-\frac{1}{p-1}\Bigl(\frac{2(p+1)}{(p-1)^2}\Bigr)^{\frac{2}{p-1}}
\quad{\rm and}\quad
H'=\Bigl(\frac{p+3}{p-1}\rho-\frac{2}{\rho}\Bigr)u'^2\leq 0\,.
\end{equation}
Hence, for solutions regular at $\rho=0$ we get
\begin{equation}
H(\rho)\leq H(0)=\frac{c^{p+1}}{p+1}-\frac{(p+1)c^2}{(p-1)^2}
\end{equation}
which implies
\begin{equation}\label{bound1} \sqrt{1-\rho^2}\,|u'(\rho)|\leq
c^{\frac{p+1}{2}}\,,
\end{equation}
and \begin{equation}\label{bound2} |u(\rho)|\leq c\quad{\rm
for}\quad c\geq (2(p+1)/(p-1)^2)^{1/(p-1)}\,. \end{equation}
Furthermore, $H(\rho)$ and thus $u(\rho)$ has a finite limit at
$\rho=1$.

On the other hand, for solutions regular at $\rho=1$ we may
integrate the function $H$ to the left. It is easy to see that
\begin{equation}
    -(H+1)^{-1}H'\leq 2/\rho,
\end{equation}
 hence $H$ stays finite on $0<\rho\leq 1$ and with it $u$ and $u'$.
 We emphasize that this argument does not exclude solutions which are
 regular at $\rho=1$ and singular at $\rho=0$. Actually, such an
 explicit singular solution
 \begin{equation}\label{uinf}
u_\infty(\rho)=b_\infty\rho^{-\frac{2}{p-1}}, \quad
b_\infty=\Biggl(\frac{2(p-3)}{(p-1)^2}\Biggr)^{\frac{1}{p-1}},
\end{equation}
will play an important role in our analysis. Note that $u_{\infty}$
corresponds to the singular static solution $\Phi=b_\infty
r^{-\frac{2}{p-1}}$ of equation (\ref{deq}).

Next, we consider the behavior of  solutions regular at $\rho=0$,
i.e. solutions satisfying the initial condition (\ref{local0}), with
large values of $c=u(0)$. To this order we rescale the variables
\begin{equation}
    \rho=\frac{x}{c^{(p-1)/2}},\quad u(\rho)=cU(x).
\end{equation}
From equation (\ref{sseq}) we get
\begin{equation}\label{sceq} \frac{d^2U}{dx^2}+\frac{2}{x}\frac{dU}{dx}+U^p=
    \frac{1}{c^{p-1}}\Bigl(x^2\frac{d^2U}{dx^2}+
  \frac{2(p+1)}{p-1}x\frac{dU}{dx}+
   \frac{2(p+1)}{(p-1)^2}U\Bigr)\,.
\end{equation}
In the limit $c\to\infty$ we obtain the simple equation
\begin{equation}\label{limeq} \frac{d^2U}{dx^2}+\frac{2}{x}\frac{dU}{dx}+U^p=0
\end{equation}
on the interval $0\leq x<\infty$ together with the boundary
condition at $x=0$
\begin{equation}\label{bdc} U(x)=1-x^2/6+O(x^4). \end{equation}
Note that the special solution $u_\infty$ is invariant under the
performed rescaling and thus $U_\infty(x)=b_\infty x^{-2/(p-1)}$
solves equation (\ref{limeq}). From our numerical analysis we know
that $U$ oscillates around $U_\infty$ suggesting the change of
variables
\begin{equation}\label{ubar}
\bar U(\tau)=x^{2/(p-1)}U(x), \quad \tau=\ln(x).
\end{equation}
Substituting (\ref{ubar}) into equation (\ref{limeq}) we get the
autonomous equation
\begin{equation}\label{limaut} \frac{d^2\bar U}{d\tau^2}+\frac{p-5}{p-1}\frac{d\bar
U}{d\tau}+\bar{U}^p
   -\frac{2(p-3)}{(p-1)^2}\bar U=0\,.
\end{equation}
Using the Lyapunov function
\begin{equation}\label{taueq} h=\frac{1}{2}\frac{d\bar
U}{d\tau}^2+\frac{\bar{U}^{p+1}}{p+1}
  -\frac{(p-3)\bar{U^2}}{(p-1)^2}
\quad{\rm with}\quad \frac{dh}{d\tau}=
  -\frac{p-5}{p-1}\Bigl(\frac{d\bar U}{d\tau}\Bigr)^2\,,
\end{equation}
we conclude that for $\tau\to\infty$ the solution $\bar U(\tau)$
tends to the fixed point $\bar U_\infty=b_\infty$ of equation
(\ref{limaut}) and $\frac{d\bar U}{d\tau}\to 0$ for $\tau\to\infty$.
Putting $\bar U=b_\infty+y$ and $\frac{d\bar U}{d\tau}=z$ we find
\begin{subequations}
\begin{eqnarray}\label{lin}
\frac{dy}{d\tau}&=&z\,,\\
\frac{dz}{d\tau}&=&-\frac{p-5}{p-1}z-\frac{2(p-3)}{p-1}y+N(y)\,,
\end{eqnarray}
\end{subequations}
where $N(y)$ denotes non-linear terms in $y$. Neglecting these
non-linear terms we get a linear system with the eigenvalues
$\lambda=(-p+5\pm i\sqrt{7p^2-22p-1})/2(p-1)$. This implies that the
solution of equation (\ref{limeq}) which is regular at $x=0$ has for
$x\to\infty$ the asymptotic form
\begin{equation}\label{asy}
U(x)=b_{\infty} x^{-\frac{2}{p-1}}\Bigl(1+A_0
x^{-\frac{p-5}{2(p-1)}}\sin(\omega\ln{x}
   +\delta_0)\Bigr)
\end{equation}
 with some constants $A_0$ and $\delta_0$ and
$\omega=\sqrt{7p^2-22p-1}/2(p-1)$.

In order to prove the existence of regular solutions for
sufficiently large $c$ we choose some $x_0$, such that for $x>x_0$
we can safely use the asymptotic expression (\ref{asy}) for the
solution of equation (\ref{limeq}) with the boundary condition
(\ref{bdc}). Furthermore, we choose $c_0$ so large that we can
neglect the r.h.s.\ of equation (\ref{sceq}) on the interval $0\leq
x<x_0$ for $c>c_0$. Integrating solutions regular at $\rho=0$  to
$\rho_0=x_0/c_0^{(p-1)/2}$ and varying the parameter $c=u(0)$  we
get a smooth curve $C_0$ in the $(u(\rho_0),u'(\rho_0))$ plane.
Using the asymptotic expression (\ref{asy}) for the part of $C_0$
corresponding to $c>c_0$ we get
\begin{equation}\label{ccurve}
    u(\rho) \approx b_{\infty} \rho^{-\frac{2}{p-1}}\Bigl(1+A_0
    c^{\frac{5-p}{4}}
\rho^{-\frac{p-5}{2(p-1)}}\sin(\omega\ln(c^{\frac{p-1}{2}}\rho)
   +\delta_0)\Bigr)\,.
\end{equation}
Keeping $\rho=\rho_0\ll 1$ fixed and letting $c\to\infty$ (and with
it $x$) we find that $C_0$ spirals down (in the negative direction)
to the limit point
\begin{equation}
P(\rho_0) =
(b_\infty\rho_0^{-2/(p-1)},-\frac{2}{p-1}b_\infty\rho_0^{-(p+1)/(p-1)}).
\end{equation}
Having that it is easy to show existence of regular solutions for
large $c$. Integrating the solutions regular at $\rho=1$ back to
$\rho_0$ we obtain a smooth curve $C_1$ parametrized by $b$. For
$0\leq b\leq b_0$ this curve goes from the point $(0,0)$ to
$(b_0,0)$ passing through $P(\rho_0)$ for $b=b_{\infty}$. Thus, the
curve $C_0$ intersects the curve $C_1$ again and again as $c$ goes
to infinity (note that the curve $C_1$ has positive rotation around
$P(\rho_0)$). This is illustrated in Fig.~1 in the case $p=7$.
Obviously each such intersection yields a solution which is regular
on the whole interval $0\leq \rho\leq 1$.
 As $c$ and thus the domain for $x$
increases the corresponding solutions $U(x)$ perform more and more
oscillations around $U_{\infty}(x)$. In fact, each further
oscillation of $\bar U(x)$ about zero leads to another cycle in the
$(u,u')$ plane at $\rho_0$. Thus, we arrive at the conclusion that
for any sufficiently large  $n$ there exist regular solutions of
equation (\ref{sseq}) on the interval $0\leq \rho\leq 1$ with $n$ of
such oscillations.
\begin{figure}[h]
\centering
\includegraphics[width=0.6\textwidth]{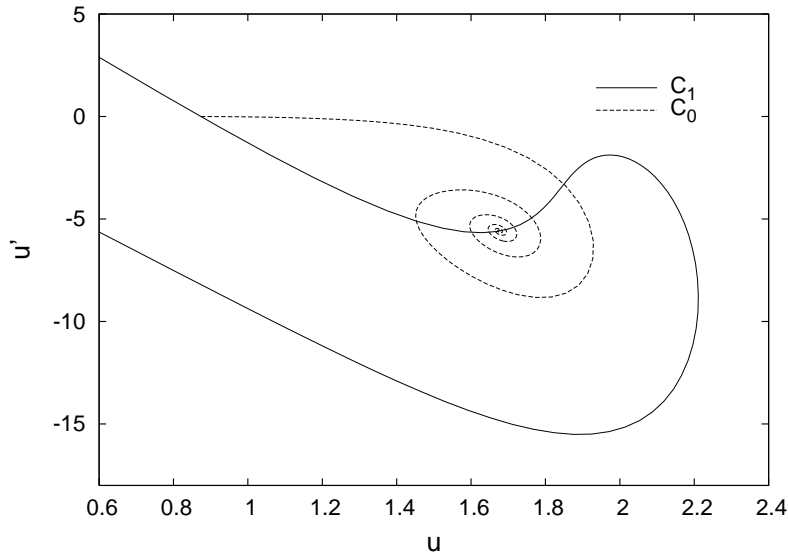}
\caption{\small{The intersection of curves $C_0$ and $C_1$ for
$\rho_0=0.1$ ($p=7$).} }\label{fig1}
\end{figure}

In the second step of our existence proof we shall show that we can
decrease $n$ and thus prove existence of regular solutions for any
positive integer $n$. The constant solution $u_0$ with $c=c_0=b_0$
has exactly one such intersection. Thus we must be able to decrease
the number of intersections to one decreasing $c$ from some large
value $c\gg 1$ corresponding to a regular solution with $n\gg 1$
intersections.
 Obviously the number of intersections can only
change at $\rho=1$, because, by the uniqueness theorem for solutions
of differential equations, a double zero of
$u(\rho,c)-u_\infty(\rho)$ is not possible at a regular point
$0<\rho<1$. On the other hand, no solution regular on the whole
interval $0\leq\rho\leq 1$ can have $b=u(1)=b_\infty$. Thus, a
change of $n$ can only happen for a value of $c$ corresponding to a
solution which is singular at $\rho=1$. We showed in Section~2 that
such solutions are parametrized by the numbers $b=u(1)$ and
$d=V(1)$, where the latter is the finite coefficient of the singular
mode. For solutions regular at $\rho=0$ the coefficient $d$ depends
smoothly on $c$. Suppose $b(c)=b_\infty$ for some $c=c_s$. Then, we
get a decrease of $n$ exactly by one  if either $b-b_\infty$ moves
from positive to negative values as $c$ decreases through $c_s$ and
$d(c_s)<0$, or $b-b_\infty$ moves from negative to positive values
as $c$ decreases through $c_s$ and $d(c_s)>0$.
This shows that as we decrease $c$ from the value corresponding to a
regular solution with $k \gg 1$ intersections to $c_0$ we
necessarily encounter regular solutions with any value of $n$
between one and $k$. This concludes the proof of Theorem~1.

 Clearly,
there could be more than one regular solution with the same number
of intersections, because $b-b_\infty$ can have zeros of even order.
However, our numerical analysis shows that there is exactly one
regular solution for any $n$. The first few solutions generated
numerically are shown in Fig.~2.
\begin{figure}[h]
\centering
\includegraphics[width=0.6\textwidth]{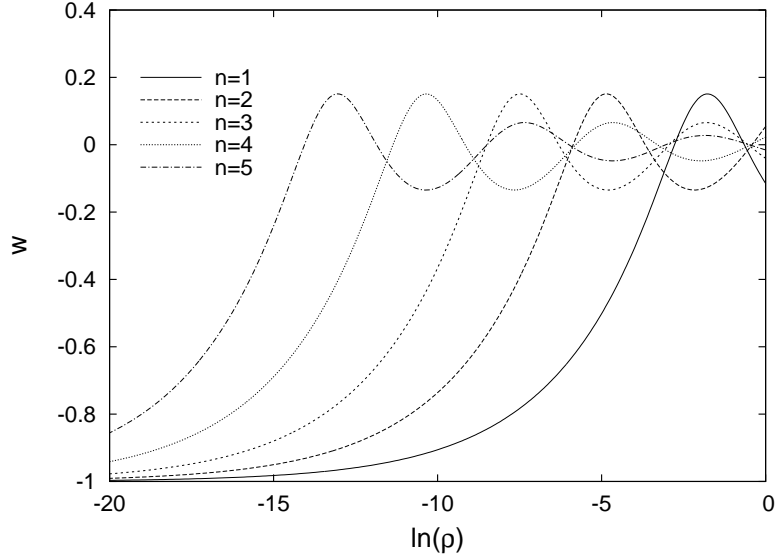}
\caption{\small{The first five self-similar solutions for $p=7$. To
show oscillations around $u_{\infty}$ we plot
$w_n=u_n/u_{\infty}-1$.} }\label{fig1}
\end{figure}
\subsection*{4. Extension beyond the past light cone}
\noindent In this section we show that the solutions constructed
above can be smoothly extended beyond $\rho=1$ to infinity.

 Consider the function
\begin{equation}\label{q}
    Q(\rho) = \frac{1}{2} (1-\rho^2) \rho^3 u'^2 + \frac{1}{2}
    \rho^2 (1-\rho^2) u u' +
    \left[\frac{3 (5-p)}{4 (p-1)}-\frac{2}{(p-1)^2}\right] \rho^3 u^2 +
    \frac{1}{p+1} \rho^3 u^{p+1}.
\end{equation}
This function was  introduced by Kavian and Weissler \cite{kw} in
their study of equation (\ref{sseq}). By straightforward (but
tedious) computation one gets
\begin{equation}\label{dq}
    Q'(\rho)=\frac{5-p}{2(p-1)} \rho^2 \left[2(3u/2+\rho
    u')^2+\frac{3p-7}{2(p-1)} u^2 + \frac{p-1}{p+1} u^{p+1}\right],
\end{equation}
hence $Q(\rho)$ is monotone decreasing for $p> 5$. Since $Q(0)=0$
for analytic solutions, it follows that $Q(\rho)\leq 0$ for $\rho>0$
and thus $u(\rho)$ is positive for $0<\rho\leq 1$. \vskip 0.2cm
\noindent \emph{Remark.} We note in passing that for $p=5$  it
follows from (\ref{dq}) that the function $Q(\rho)$ is the first
integral which implies that $u_0$ is the only regular self-similar
solution. \vskip 0.2cm
 It is convenient to rewrite equation
(\ref{sseq}) in a self-adjoint form
\begin{equation}\label{sadjoint}
    (\rho^2(1-\rho^2)^{\alpha} u')'=\rho^2 (1-\rho^2)^{\alpha-1} u
    (b_0^{p-1}-u^{p-1}).
\end{equation}
\textbf{Proposition~1.} If $u(\rho)$ is a regular solution of
equation (\ref{sadjoint}) defined for $0\leq \rho\leq 1$ and
$u(0)>b_0$, then $u(1)<b_0$. \vskip 0.1cm \noindent \emph{Proof:}
First, note that if $u$ is analytic at $\rho=1$ and $u(1)=b_0$, then
$u(\rho)\equiv b_0$, which contradicts the assumption that
$u(0)>b_0$. Thus, the case $u(1)=b_0$ is impossible. Next, suppose
that $u(1)>b_0$. If $u(\rho)\geq b_0$ on the whole interval $0\leq
\rho\leq 1$ then the r.h.s of equation (\ref{sadjoint}) is negative
and hence by integrating it we get that $u'(1)=\infty$. Thus,
$u(\rho)<b_0$ for some $0<\rho<1$ and therefore $u$ must have a
minimum value at some $\rho_0$, i.e., $u'(\rho_0)=0$ and
$0<u(\rho_0)<b_0$. If so, between $\rho_0$ and $1$ there is a
$\rho_1$ where $u(\rho_1)=b_0$ and $u(\rho)<b_0$ for $\rho_0\leq
\rho<\rho_1$. We now show that this is not possible. To this order
we integrate equation (\ref{sadjoint}) over the interval $\rho_0\leq
s<\rho< \rho_1$. Since $u'(\rho_0)=0$ we get
\begin{equation}\label{p1}
\rho^2(1-\rho^2)^{\alpha} u'(\rho)=\int_{\rho_0}^{\rho} s^2
(1-s^2)^{\alpha-1} u
    (b_0^{p-1}-u^{p-1}) ds.
\end{equation}
We have $u(s)<u(\rho)<b_0$, hence
$
u(b_0^{p-1}-u^{p-1})=u (b_0^{p-2}+b_0^{p-3} u + ...+ u^{p-2})
(b_0-u) < (p-1) b_0^{p-1}(b_0-u)$. Thus,
\begin{equation}\label{cal}
    \int_{\rho_0}^{\rho} s^2
(1-s^2)^{\alpha-1} u
    (b_0^{p-1}-u^{p-1}) ds <(p-1) b_0^{p-1}(b_0-u(\rho)) \int_{\rho_0}^{\rho} s^2
(1-s^2)^{\alpha-1} ds = f(\rho) (b_0-u(\rho))
\end{equation}
for some bounded continuous function $f(\rho)$. It follows from
(\ref{p1}) and (\ref{cal}) that
\begin{equation}\label{p2}
    \frac{u'(\rho)}{b_0-u(\rho)}
    <\frac{f(\rho)}{\rho^2(1-\rho^2)^{\alpha}}
    =g(\rho),
\end{equation}
where $g(\rho)$ is a bounded continuous function on $\rho_0<\rho\leq
\rho_1<1$. Integrating (\ref{p2}) from $\rho_0$ to $\rho$, we get
$\ln(b_0-u(\rho_0))-\ln(b_0-u(\rho))<\int_{\rho_0}^{\rho} g(s) ds$,
hence $\ln(b_0-u(\rho))$ is bounded for $\rho\leq \rho_1$. This
contradicts the assumption that $u(\rho_1)=b_0$ and concludes the
proof. \vskip 0.1cm  Next, we consider the behavior of solutions for
$\rho>1$. \vskip 0.1cm \noindent \textbf{Proposition~2.} If
$u(\rho)$ is a solution to equation (\ref{sadjoint}) that is
analytic at $\rho=1$ and $u(1)=b<b_0$, then $u(\rho)$ is defined for
all $\rho>1$ and $0<u(\rho)<b_0$. \vskip 0.1cm \noindent
\emph{Proof:} Integrating equation (\ref{sadjoint}) from $1$ to some
$\rho>1$ we get
\begin{equation}\label{con}
    \rho^2 (\rho^2-1)^{\alpha} u'(\rho) =\int_1^{\rho} s^2 (1-s^2)^{\alpha} u
    (b_0^{p-1}-u^{p-1}) ds,
\end{equation}
hence if $u$ is bounded, so is $u'$ and therefore the solution
cannot become singular. Thus, it suffices to show that
$0<u(\rho)<b$. From equation (\ref{sadjoint}) we have that
$u'(\rho)=0$ and $0<u(\rho)<b_0$ for $\rho>1$ implies that
$u''(\rho)<0$, thus $u$ is monotone decreasing. To show that $u$
cannot become negative we introduce the function $f(\rho)=\rho
u'(\rho)+(\alpha+1) u(\rho)/2$. Evaluating $f'$ at $f=0$ we obtain
\begin{equation}\label{fprim}
    f'\Bigr\rvert_{f=0} = \frac{1}{2} (1-\alpha^2) u +
    \frac{2\rho^2}{\rho^2-1} u^p >0,
\end{equation}
which, together with  $f(1)=b^p/2\alpha>0$, implies that $f(\rho)>0$
for $\rho>1$. Since $u'(\rho)<0$ (as noted above), we conclude that
$u>0$, as claimed.

\subsection*{5. Large $n$ asymptotics}
\noindent Our numerical analysis shows that the shooting parameters
$b_n$ and $c_n$ exhibit remarkable scaling properties in the large
$n$ limit.
 In this section we
explain this phenomenon.

In order to derive the large $n$ asymptotic behavior of  solutions
which are regular at $\rho=1$ we introduce a new variable $w$
defined by $u=u_{\infty} (1+w)$. Using  equations (\ref{sseq}) and
(\ref{local1}) we find that $w$ satisfies the equation
\begin{equation}\label{weq}
    \rho^2 (1-\rho^2) w'' +\left(\frac{2(p-3)}{p-1} \rho -
    2\rho^3\right) w' -\frac{2(p-3)}{(p-1)^2} (1+w)
    ((1+w)^{p-1}-1)=0,
\end{equation}
with the boundary condition at $\rho=1$
\begin{equation}\label{weg_bc}
    w(\rho)=\frac{b-b_{\infty}}{b_{\infty}} +
    \frac{p-1}{4}\frac{b}{b_{\infty}} (b^{p-1}-b_{\infty}^{p-1})
    (\rho-1) + O((\rho-1)^2).
\end{equation}
For any $0<\rho_0\ll 1$ and $b$ sufficiently close to $b_{\infty}$,
the solution of equation (\ref{weq}) satisfying the initial
condition (\ref{weg_bc})  will stay arbitrarily close to $w=0$ on
the interval $\rho_0\leq\rho\leq 1$, hence it can be approximated by
$w=\frac{b-b_{\infty}}{b_{\infty}} w_L$, where $w_L$ satisfies the
linearized equation
\begin{equation}\label{weqlin}
    \rho^2 (1-\rho^2) w_L'' +\left(\frac{2(p-3)}{p-1} \rho -
    2\rho^3\right) w_L' -\frac{2(p-3)}{(p-1)} w_L=0,
\end{equation}
with the boundary condition $w_L(1)=1, w_L'(1)=(p-3)/2$. We could
express $w_L$ in terms of a hypergeometric function but for our
purposes it is sufficient to have its asymptotics for small values
of $\rho$
\begin{equation}\label{wl}
    w_L\approx \rho^{\frac{5-p}{2(p-1)}}
    A_1\sin(\omega\ln{\rho}+\delta_1),\qquad
    \omega=\sqrt{7p^2-22p-1}/2(p-1)\,,
\end{equation}
where $A_1$ and $\delta_1$ are constants. Thus, the solution regular
at $\rho=1$ with $b$ near $b_{\infty}$ has the following form for
small $\rho$
\begin{equation}\label{uas}
    u(\rho)\approx u_{\infty} \left(1+\frac{b-b_{\infty}}{b_{\infty}}
    \rho^{\frac{5-p}{2(p-1)}}
    A_1\sin(\omega\ln{\rho}+\delta_1)\right).
\end{equation}
The key point is that the regions of validity of the approximations
(\ref{uas}) and (\ref{ccurve}) overlap so we can match them.
Matching the amplitudes of solutions (\ref{ccurve}) and (\ref{uas})
at $\rho_0$ we get the condition
\begin{equation}\label{amatch}
    A_0 c^{\frac{p-5}{4}} = \pm A_1 \frac{b-b_{\infty}}{b_{\infty}}.
\end{equation}
On the other hand, using (\ref{ccurve}) and the fact that  the
phases of two adjacent solutions $u_{n}(\rho)$ and $u_{n+1}(\rho)$
differ by $\pi$ we obtain
\begin{equation}\label{pmatch}
    \frac{p-1}{2} \omega \ln\left(\frac{c_{n+1}}{c_n}\right) = \pi,
\end{equation}
hence
\begin{equation}\label{scal1}
    \frac{c_{n+1}}{c_n} \approx e^{\frac{2\pi}{(p-1)\omega}}\,,
\end{equation}
which, using (\ref{amatch}), implies that
\begin{equation}\label{scal2}
    \frac{b_{n+1}-b_{\infty}}{b_{\infty}-b_n} \approx
    e^{-\frac{(p-5)\pi}{2(p-1)\omega}}\,.
\end{equation}
The numerical verification of these scaling properties is shown in
Table~1 in the case $p=7$.

\vskip 0.2cm \setlength{\tabcolsep}{1.5em}
\begin{table}[h]
\begin{center}
\begin{tabular}
{|c|c|c|c|c|}
  \hline
  $n$ & $c_n$ & $b_n$ & $\Delta c_n$ & $\Delta b_n$ \\
  \hline
  1 & 2.054390385 & 0.688698572 & 2.8018 & 0.4713 \\
  2 & 5.756037116 & 0.820493408 & 2.4090 & 0.7428 \\
  3 & 13.86655615 & 0.746908360 & 2.5899 & 0.5670 \\
  4 & 35.91343330 & 0.796055093 & 2.4577 & 0.6752 \\
  5 & 88.26661166 & 0.766263419 & 2.5326 & 0.6059 \\
  6 & 223.5507381 & 0.785548198 & 2.4823 & 0.6493 \\
  7 & 554.9215495 & 0.773546658 & 2.5128 & 0.6217 \\
  8 & 1394.439242 & 0.781209584 & 2.4930 & 0.6391 \\
  9 & 3476.402010 & 0.776393971 & 2.5053 & 0.6281 \\
  10 & 8709.676250 & 0.779451184 & 2.4974 & 0.6351 \\
  11 & 21752.40861 & 0.777522645 & 2.5012 & 0.6307 \\
  12 & 54434.14714 & 0.778744142 & 2.4993 & 0.6334 \\
  13 & 136047.6759 & 0.777972446 & 2.5000 & 0.6317 \\
  14 & 340293.1022 & 0.778460765 & 2.5008 & 0.6328 \\
  15 & 850746.1358 & 0.778152079 & 2.5003 & 0.6321 \\ \hline
  $\infty$ & $\infty$ & 0.778271716 & 2.5005 & 0.6324 \\
  \hline
\end{tabular}
  \caption{Shooting parameters of solutions $u_n$ and their quotients $\Delta c_n=c_{n+1}/c_n$ and
  $\Delta b_n=(b_{n+1}-b_{\infty})/(b_{\infty}-b_n)$ for $p=7$. The last row for $n=\infty$
  corresponds to the analytic results (45) and (46).}
  \end{center}
\end{table}
\vskip 0.2cm \noindent \emph{Note added in proof.} Peter
Breitenlohner (personal communication) pointed out that a minor
modification of our proof applies also to even values of $p$, hence
the conclusion of Theorem~1 holds for all integers $p\geq 6$. \vskip
0.2cm

\subsection*{Acknowledgments} We are grateful to Peter Breitelohner
for the careful reading of the manuscript and comments which helped
us to correct and improve it.  We thank the referee for pointing out
the reference \cite{l} where similar ideas to that of section~3 were
used before to prove existence of self-similar solutions of a
nonlinear heat equation.  PB and AW thank the Mathematisches
Forschungsinstitut in Oberwolfach for supporting this project under
the "Research in Pairs" program. The work of PB was supported in
part by the Polish Ministry of Science grant no. 1PO3B01229.

\section{Appendix}
\noindent We present here an alternative  version of the proof of
Theorem~1. Throughout the appendix the range of $\rho$ is $0\leq
\rho \leq 1$. \vskip 0.2cm  We know from Section~2 and equation (16)
that given any $c$ there is a unique solution $u(\rho,c)$ of
equation (5) satisfying $u(0,c)=c$ defined for all $0\leq \rho\leq
1$. For $\rho>0$ we define $v(\rho,c)=u(\rho,c)/u_{\infty}(\rho)$
and $w(\rho,c)=v(\rho,c)-1$. We  define also the function
\begin{equation}\label{hv}
    H_v=\frac{1}{2} \rho^2(1-\rho^2)
    v'^2-\alpha(1-\alpha)\left(\frac{v^2}{2}-\frac{v^{p+1}}{p+1}\right)\,.
\end{equation}
Using equation (5) we get
\begin{equation}
    H_v'=\frac{5-p}{p-1} \rho v'^2\,,
\end{equation}
hence $H_v$ is decreasing along solutions. Furthermore, $H_v$ has a
minimum equal to $-(p-3)/(p^2-1)$ for $v'=0$ and $v=1$. \vskip 0.1cm
\noindent
\textbf{Lemma~1} Let $d=(p-1)/p$. If $c>b_0/d$, then there is an
$\rho_1=\rho_1(c)<(b_{\infty}/dc)^{1/\alpha}$ such that
$w(\rho_1,c)=0$ and $0<\rho_1 w'(\rho_1,c)<\alpha$. \vskip 0.2cm
\noindent \emph{Proof:}
It follows from equation (32) that $u(\rho,c)$ is a decreasing
function of $\rho$ for $c>b_0$ as long as $u>b_0$. Integrating
equation (32) and noting that $0<u(\rho,c)\leq c$  we get
\begin{equation}\label{32int}
    \rho^2(1-\rho^2)^{\alpha} u'(\rho)=-\int_0^{\rho}
    s^2(1-s^2)^{\alpha-1}[u(s)^p-\alpha(1+\alpha) u(s)] ds > -\int_0^{\rho}
    s^2(1-s^2)^{\alpha-1}u(s)^p ds >  -\int_0^{\rho}
    s^2(1-s^2)^{\alpha-1}c^p ds\,.
\end{equation}
Thus
\begin{equation}
    \rho^2(1-\rho^2)^{\alpha} u'(\rho) >  -\int_0^{\rho}
    \frac{s^2(1-s^2)^{\alpha-1}}{(1-\rho^2)^{\alpha-1}} c^p ds
    >-\int_0^{\rho} s^2 c^p ds =-\frac{\rho^3}{3} c^p\,,
\end{equation}
hence
\begin{equation}\label{up}
    u'(\rho)>-\frac{\rho}{3(1-\rho^2)} c^p\,.
\end{equation}
Integrating (\ref{up}) from $0$ to $x$ where $u(x,c)=d c$ gives
\begin{equation}
    -c(1-d)>\frac{c^p}{6} \ln(1-x^2),\qquad \mbox{so} \qquad
    x^2>1-\exp\left(-\frac{6}{pc^{p-1}}\right) > \frac{6}{p c^{p-1}}\,.
\end{equation}
Let $y$ be a point where $u_{\infty}(y)=b_{\infty} y^{-\alpha}=dc$,
that is, $y=(b_{\infty}/dc)^{1/\alpha}$. We have
\begin{equation}
    \left(\frac{b_{\infty}}{d}\right)^{2/\alpha} =
    \left(\frac{b_{\infty}}{d}\right)^{p-1}=\left(\frac{p}{p-1}\right)^{p-1}\frac{2(p-3)}{(p-1)^2}<\frac{6}{p}\,,
\end{equation}
hence $x>y$, and therefore $u(y,c)>u(x,c)=dc=u_{\infty}(y)$, or
$w(y,c)>0$. Since $w(\rho,c)<0$ for small $\rho$, we conclude that
there is a $\rho$ such that $w(\rho,c)=0$ and, if $\rho_1$ is the
smallest such $\rho$, then $w'(\rho_1,c)>0$. Note that $\rho_1<y$.
This conludes the proof of Lemma~1. \vskip 0.2cm \noindent
\textbf{Corollary.} $w(\rho)>-2\alpha$ for $\rho>\rho_1$. \vskip
0.1cm \noindent \emph{Proof:} Using the fact that the function
$H_v(\rho)$ is monotone decreasing we obtain
\begin{equation}\label{cor}
    H_v(\rho)\leq H_v(\rho_1)\leq
    \frac{\alpha^2}{2}-\alpha(1-\alpha)\left(\frac{1}{2}-\frac{1}{p+1}\right)\,,
\end{equation}
for $\rho\geq \rho_1$. On the other hand, for  $\rho\geq \rho_1$
\begin{equation}\label{cor2}
H_v(\rho)=\frac{1}{2} \rho^2(1-\rho^2)
    v'^2-\alpha(1-\alpha)\left(\frac{v^2}{2}-\frac{v^{p+1}}{p-1}\right)\geq
    -\alpha(1-\alpha)\left(\frac{v^2}{2}-\frac{v^{p+1}}{p-1}\right)\,.
\end{equation}
Combining equations (\ref{cor}) and (\ref{cor2}) we have
\begin{equation}
f(v)\equiv
\frac{\alpha^2}{2}-\alpha(1-\alpha)\left(\frac{1}{2}-\frac{1}{p+1}\right)
    +\alpha(1-\alpha)\left(\frac{v^2}{2}-\frac{v^{p+1}}{p-1}\right)
    \geq 0\,,
    \end{equation}
    for $\rho\geq \rho_1$. Note that $f'(v)>0$ for $v<1$. A straightforward calculation yields
    $f(1-2\alpha)<0$ if $p\geq 7$, hence $v>1-2\alpha$ (or
    $w>-2\alpha$) for $\rho\geq \rho_1$, as claimed.
\vskip 0.2cm Let us define
\begin{equation}\label{R}
    R(\rho,c)=\sqrt{w(\rho,c)^2+w'(\rho,c)^2}\,, \qquad 0<\rho\leq
    1\,.
\end{equation}
By uniqueness of solutions of ODEs, a solution starting at
$(u,u')=(c,0)$ at $\rho=0$ cannot have
$(u(\rho,c),u'(\rho,c)=(u_{\infty}(\rho),u'_{\infty}(\rho))$, which
implies that $R(\rho,c)>0$. Thus, we may also define
\begin{equation}\label{theta}
    \Theta(\rho,c)=\arctan\left(\frac{\rho
    w'(\rho,c)}{w(\rho,c)}\right)\qquad \mbox{up to a multiple of
    $\pi$}\,.
\end{equation}
Since the region $\{(\rho,c)| 0<\rho\leq 1,c\geq 0\}$ is simply
connected we may unambiguously define a function $\Theta(\rho,c)$
once we specify its value at any point in the domain. Since
$w(\rho,0)\equiv -1, w'(\rho,0)\equiv 0$, we set
$\Theta(1/2,0)=\pi$. Note that $\lim_{\rho\rightarrow 0+}
w(\rho,c)=-1$ and $\lim_{\rho\rightarrow 0+} \rho w'(\rho,c)=0$,
hence $\lim_{\rho\rightarrow 0+} \Theta(\rho,c)=\pi$. \vskip 0.2cm
\noindent \textbf{Lemma~2.} Let $\rho_2\ll 1$. We have
$\lim_{c\rightarrow \infty}\Theta(\rho_2,c)= -\infty$. \vskip 0.2cm
\noindent \emph{Proof:} Consider the equation satisfied by the
function $\Theta(\rho)=\Theta(\rho,c)$:
\begin{equation}\label{Thetaeq}
    \Theta'(\rho)=-\frac{1}{\rho}\left[\sin^2{\Theta}+\frac{(p-5)-\rho^2(p-1)}{(p-1)(1-\rho^2)}\;\sin{\Theta}\cos{\Theta}
    +\frac{2(p-3)v}{(p-1)^2(1-\rho^2)}
    (1+v+...+v^{p-2})\cos^2{\Theta}\,.
    \right]
\end{equation}
We will show that the quantity in brackets in bounded from below by
a positive constant $\eta$ for $\rho\geq \rho_1$. It follows that
\begin{equation}
\lim_{c\rightarrow \infty}\Theta(\rho_2,c)=\lim_{c\rightarrow
\infty}\left(\Theta(\rho_1(c),c)+\int_{\rho_1(c)}^{\rho_2}
\Theta'(\rho) d\rho\right)\leq \frac{\pi}{2}+\lim_{c\rightarrow
\infty} \int_{\rho_1(c)}^{\rho_2} \left(-\frac{\eta}{\rho}\right)=
\frac{\pi}{2}-\eta \lim_{c\rightarrow
\infty}\ln\frac{\rho_2}{\rho_1(c)}=-\infty\,,
\end{equation}
and the lemma is proven.

To get a bound on the quantity in brackets we regard it as a
(formal) quadratic form in $\sin{\Theta}$ and $\cos{\Theta}$:
$[\,]=A \sin^2{\Theta}+B\sin{\Theta}\cos{\Theta}+C\cos^2{\Theta}$,
where
\begin{equation}
   A=1,\qquad  B=\frac{(p-5)-\rho^2(p-1)}{(p-1)(1-\rho^2)},\qquad
    C=\frac{2(p-3)v}{(p-1)^2(1-\rho^2)}
    (1+v+...+v^{p-2})\,.
\end{equation}
Since $\rho<\rho_2$, we may replace these coefficients by
\begin{equation}
   A=1,\qquad  \tilde B=\frac{(p-5)}{(p-1)},\qquad
    \tilde C=\frac{2(p-3)v}{(p-1)^2}
    (1+v+...+v^{p-2})\,.
\end{equation}
We need to show that the discriminant $\Delta=\tilde B^2-4A\tilde C$
is negative for $v>1-2\alpha$. We have
\begin{equation}\label{Delta}
    \Delta =\frac{(p-5)^2}{(p-1)^2}-\frac{8(p-3)}{p-1)^2}
    (1+v+...+v^{p-2})\,.
\end{equation}
Clearly, $(p-1)^2\Delta$ is a decreasing function of $v$. When
$v=1-2\alpha$ we get
\begin{equation}
    (p-1)^2\Delta =
    (2p^2-8p+6)\left(\frac{p-5}{p-1}\right)^p-p^2+6p-5\,,
\end{equation}
which is negative for $p\geq 7$, hence $\Delta<0$ which concludes
the proof of Lemma~2. \vskip 0.2cm Next, we consider solutions
$U(\rho,b)$ which start from initial value (10). As we showed in
Section~3 these solutions are defined for $0<\rho\leq 1$, hence we
may define
\begin{equation}\label{Rb}
    \tilde R(\rho,b)=\sqrt{W(\rho,b)^2+W'(\rho,b)^2}\,, \qquad 0<\rho\leq
    1\,,
\end{equation}
where $W(\rho,b)=U(\rho,b)/u_{\infty}(\rho)-1$. By uniqueness of
solutions of ODEs, a solution starting at $\rho=1$ cannot have
$(U(\rho,b),U'(\rho,b)=(u_{\infty}(\rho),u'_{\infty}(\rho))$, unless
$b=b_{\infty}$, which implies that $\tilde R(\rho,b)>0$ if $b\neq
b_{\infty}$. Thus, we may also define
\begin{equation}\label{theta}
    \tilde \Theta(\rho,b)=\arctan\left(\frac{\rho
    W'(\rho,b)}{W(\rho,b)}\right)\qquad \mbox{up to a multiple of
    $\pi$}\,.
\end{equation}
Since the region $\{(\rho,b)| 0<\rho\leq 1,0\leq b\leq b_{\infty}\}$
is simply connected we may unambiguously define a function $\tilde
\Theta(\rho,b)$ once we specify its value at any point in the
domain. We set $\tilde \Theta(1,0)=\pi$.

We define a map
\begin{equation}
    \Phi: R_{+}=\{c|c\geq 0\} \rightarrow R^2_{+}=\{(x,y)\in
    R^2|y>0\}\,,\qquad \Phi(c)=(\Theta(\rho_0,c),R(\rho_0,c))\,,
\end{equation}
and a map
\begin{equation}
    \Psi_n: t\in(-\infty,b_{\infty}) \rightarrow R^2_{+}\,,\qquad \Psi_n(t)=\begin{cases} (\tilde \Theta(\rho_0,b)-2n\pi,\tilde R(\rho_0,b)) &
\text{if $0\leq
    t<b_{\infty}\,\,$ (segment 1)},
\\
(\pi-2n\pi,1-t) &\text{if $t<0\,\,$ (segment 2)}.
\end{cases}
\end{equation}
Note that if $\Psi_n(b)=\Phi(c)$ for some $b\in(0,b_{\infty})$ and
some $c$, then $\Phi(c)$ cannot lie on  segment 2 because
$w(\rho_0,c)\geq -1$ and $y\cos{x}=-y<-1$ on segment 2. Thus, if
$\Psi_n(b)=\Phi(c)$ for some $b\in(0,b_{\infty})$ and some $c$, then
the functions $u(\rho,c)$ and $U(\rho,b)$ and their derivatives
match at $\rho_0$, hence  we have a solution defined on the whole
interval $0\leq \rho\leq 1$ in the $n^{th}$ nodal class. The integer
part of $(\tilde\Theta(1,b)-\Theta(0,c))/\pi$ counts the number of
zeros of $w$. In particular, $\Psi_0(0)=\Phi(0)$ since
$u(\rho,0)\equiv U(\rho,0)\equiv 0$ is a solution in the $0^{th}$
nodal class. \vskip 0.2cm \noindent \emph{Remark:} The choices of
$\Theta(1/2,0)=\pi$ and $\tilde\Theta(1,0)=\pi$ were made so that
the indexing of solutions is correct.
\newpage
\begin{figure}[h]
\centering
\includegraphics[width=0.8\textwidth]{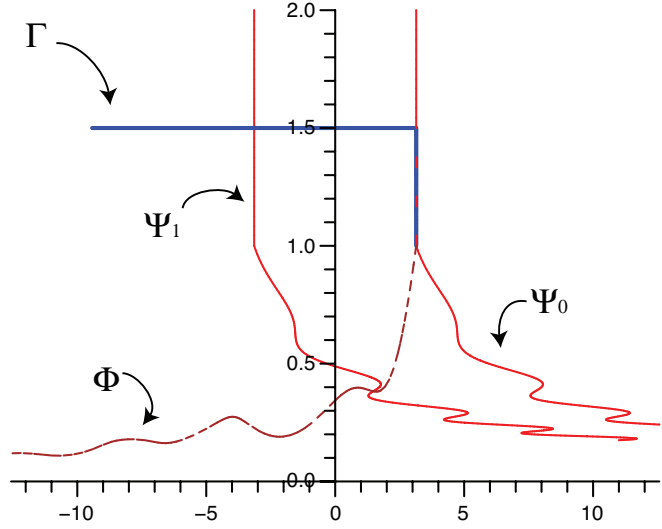}
\caption{\small{Illustration to the proof of
Proposition~3.}}\label{fig5}
\end{figure}
Now we are ready to prove: \vskip 0.2cm \noindent
\textbf{Proposition~3.}
 For any positive integer $n$ there is
$B_n<b_{\infty}$ such that the solution $U(\rho,B_n)$ is in the
$n^{th}$ nodal class. \vskip 0.2cm \noindent \emph{Proof:}  Since
$b_{\infty}<b_0$, there is an $\rho_0<1$ such that $U(\rho_0,b)<b_0$
for all $b\in[0,b_{\infty}]$. Note that $U'(1,b)<0$ for $b\leq
b_{\infty}$ and thus $U'(\rho,b)<0$ on $\rho_0<\rho<1$ because
$U''(\rho,b)>0$ if $U'(\rho,b)=0$ and $U(\rho,b)<b_0$. It follows
that $U(\rho,b)\geq 0$ and hence $W(\rho,b)\geq -1$ for
$\rho_0<\rho\leq 1$ and $0\leq b\leq b_{\infty}$. From this conclude
that the curve $\Psi_n(b)$ is simple, i.e., $\Psi_n(b)=\Psi_n(\bar
b)$ implies that $b=\bar b$. For nonnegative $b$ and $\bar b$ this
follows from the uniqueness of solutions of ODEs, while for $\bar
b<0\leq b$, $\Psi_n(\bar b)=\Psi_n(b)$ is impossible because
$W(\rho_0,b)>-1$ by choice of $\rho_0$ and $y\cos{x}=-y<-1$ on
segment 2 of the curve.
Since $\lim_{b\rightarrow b_{\infty}}\tilde R(\rho_0,b)=\tilde
R(\rho_0,b_{\infty})=0$ and $\lim_{b\rightarrow-\infty}\tilde
R(\rho_0,b)=\infty$, it follows from the Jordan curve theorem that
the curve $\Psi_n$ separates the half plane $R_+^2$ into two open
regions: $R_+^2\backslash image(\Psi_n)=A\cup B$. Moreover, points
$p$ and $q$ are in the different components $A$ and $B$ iff there is
a curve from $p$ to $q$ crossing $image(\Psi_n)$ transversally at
exactly one point.

Note that $\tilde \Theta(1,b)>\pi/2$ if $b<b_{\infty}$ and $\tilde
\Theta'(\rho,b)>0$ when $\tilde \Theta(\rho,b)=\pi/2$, hence $\tilde
\Theta(\rho,b)\geq \pi/2$ for $b<b_\infty$. The region
\begin{equation}
    Z=\{(x,y)\in R_+^2| x<\pi/2-2n\pi\}
\end{equation}
is connected and does not meet $image(\Psi_n)$ because $\tilde
\Theta(\rho,b)\geq \pi/2$, hence $Z$ must be contained in $A$ or
$B$; suppose that $Z\subset A$. The curve $\Phi$ meets $A$ because
$\lim_{c\rightarrow \infty}\Theta(\rho,c)= -\infty$ by Lemma~2 and
hence the curve $\Phi$ meets $Z\subset A$. Note that
$image(\Psi_n)\cap image(\Psi_m)=\emptyset$ if $n\neq m$ by
uniqueness of solutions of ODEs, hence $\Phi(0)\notin image(\Psi_n)$
if $n>0$.

We will now show that $\Phi(0)\in B$. That will complete the proof
because the curve $\Phi$ is connected, hence it must cross
$image(\Psi_n)$, which means that there is a solution in the
$n^{th}$ nodal class. To show this we construct a curve in $R_+^2$
from $\Phi(0)$ to $Z$ that crosses $image(\Psi_n)$ transversally at
exactly one point. We denote by $M$  an upper bound of $\tilde
R(\rho_0,b)$ on $0\leq b\leq b_{\infty}$ and define a curve
\begin{equation}
    \Gamma(t)=\begin{cases} (\pi,1+t) & \text{for $0\leq t\leq M$
    (segment 1)},\\
    (\pi+M-t,M+1) & \text{for $M\leq t\leq (2n+1)\pi$ (segment 2)}.
    \end{cases}
\end{equation}
Note that segment 1 of $\Gamma$ cannot meet segment 1 of $\Psi_n$
because segment 1 of $\Psi_n$ has $w>-1$ and segment 1 of $\Gamma$
has $y\cos{x}<-1$. Also, segment 1 of $\Gamma$ cannot meet segment 2
of $\Psi_n$ because they have different first coordinates
$\pi-2\pi<\pi$. Finally, segment 2 of $\Gamma$ cannot meet segment 1
of $\Psi_n$ because $\tilde R(\rho_0,b)\leq M<M=1$. Thus, the only
intersection point of $\Gamma$ and $\Psi_n$ is at $(\pi-2n\pi,M+1)$.
This completes the proof of Proposition~3.

\vskip 0.2cm \noindent \textbf{Proposition~4.}
 For any positive integer $n$ there is
$\tilde B_n>b_{\infty}$ such that the solution $U(\rho,b_n)$ is in
the $n^{th}$ nodal class. \vskip 0.2cm \noindent \emph{Proof:} The
proof proceeds along the same lines as the proof of Proposition~3
except for some technical modifications which we leave to the
reader. \vskip 0.2cm \noindent  Theorem~1 is the immediate
consequence of Propositions~3 and~4 if we set $b_{2n-1}=B_n$ and
$b_{2n}=\tilde B_n$.


\begin{thebibliography}{10}

\bibitem{k} J. Keller, Comm. Pure Appl. Math. \textbf{10}, 523
(1957).

\bibitem{j} F. John, Manuscripta Math. \textbf{28}, 235 (1979).

\bibitem{g} R. Glassey, Math. Z. \textbf{178}, 233 (1981).

\bibitem{dw} {\sl Dispersive PDE Wiki}, http://tosio.math.toronto.edu/wiki

\bibitem{kl} S. Kichenassamy and W. Littman, Commun. PDE
\textbf{18}, 1869 (1993).

\bibitem{bct} P. Bizo\'n, T. Chmaj, and Z. Tabor, Nonlinearity
\textbf{17}, 2187 (2004).

\bibitem{kw} O. Kavian and F. B. Weissler, Commun. PDE \textbf{15},
1421 (1990).

\bibitem{bfm}
P.~Breitenlohner, P.~Forg\'acs and D.~Maison, Comm. Math. Phys.
\textbf{163}, 141 (1994).

\bibitem{H}
P.~Hartman, {\sl Ordinary Differential Equations\/}, Boston:
Birkh\"auser, 1982.

\bibitem{l} L. A. Lepin, Differ. Equ. \textbf{24}, 799 (1988).

\end{thebibliography}
\end{document}